\documentclass{article}
\usepackage{amsfonts,amssymb,latexsym,amsmath}
	\newcounter{num}[section]
	\newcommand{\Num}{\refstepcounter{num}%
		\textbf{\arabic{section}.\arabic{num}}}
	
	\newcommand{\Theorem}{\textbf{Theorem~}}
	\newcommand{\Proof}{\textbf{Proof}}

	\newcommand{\Lemma}{\textbf{Lemma~}}
	
	\newcommand{\Remark}{\textbf{Remark~}}

	\newcommand{\ux}{{\mathfrak u}}
	\newcommand{\nx}{{\mathfrak n}}
	\newcommand{\mx}{{\mathfrak m}}
	\newcommand{\px}{{\mathfrak p}}
	\newcommand{\gx}{{\mathfrak g}}
	\newcommand{\hx}{{\mathfrak h}}
	\newcommand{\bx}{{\mathfrak b}}
	\newcommand{\Sx}{{\mathfrak S}}
	\newcommand{\De}{{\Delta}}
	\newcommand{\Dp}{{\Delta}^+}
	
	\newcommand{\al}{{\alpha}}
	\newcommand{\Ac}{{\cal A}}
	\newcommand{\Fc}{{\cal F}}
	\newcommand{\Sc}{{\cal S}}
	
	\newcommand{\Res}{{\mathrm{Res}}}
	\newcommand{\ad}{{\mathrm{ad}}}
	
	\newcommand{\wt}{{\mathrm{wt}}}

	\renewcommand{\leq}{\leqslant}
	\renewcommand{\geq}{\geqslant}
	
\begin{document}
		\Large
	\title{ $U$-projectors and   fields of $U$-invariants}
	\author{K.A.Vyatkina \and A.N.Panov
		\thanks{The paper is supported by the RFBR grants 16-01-00154-a and 14-01-97017-Volga-region-a}}
	
	\date{}
	\maketitle
	\begin{abstract}
		We present the general construction of the $U$-projector  (the homomorphism of the algebra into its field of  $U$-invariants identical on the subalgebra of $U$-invariants). It is shown how to apply   $U$-projector  to find the systems of  free generators of the fields of $U$-invariants for representations of reductive groups.
	\end{abstract}
	
	\section{Introduction}
	
	It is well known that the fields of invariants  of the triangular transformation groups of affine varieties are rational   (see. \cite{ M, Pu, V}).
	In this paper, for the special instances, we are going to present  the  systems of free generators.
	
	For this purpose we introduce the notion of $U$-projector  (it is a homomorphism of the algebra  $K[X]$ to the field of invariants $K(X)^U$ identical on  $K[X]^U$).  In the section  2,  we verify   existence of  $U$-projector and  present its general construction. Furthermore,  we show that applying the $U$-projector to the suitable system of functions  $b_1,\ldots, b_m$ one can obtain the system of functions   $P(b_1), \ldots, P(b_m)$ that freely generate  the field of  $U$-invariants.
	
	Notice that the general $U$-projector construction method  is realized by the induction method on the length of  ideals chain.
	Since the length might by rather great,  this method don't provide the exact formula for $U$-projector.
	However, a priory the   $U$-projector is not unique. In the
	next sections  3,4,5, we return to the problem of  $U$-projector construction in the special cases. Our goal is to improve the above $U$-projector construction  to make it more precise. We also plan to choose the system  of functions  $b_1,\ldots, b_m$ such that  the system   $P(b_1), \ldots, P(b_m)$ freely generate the field of  $U$-invariants. The main results are formulated in the theorems  \ref{acac},~ \ref{gen},~ \ref{atwo}, ~\ref{athree},~\ref{afour}. Notice that in the paper \cite{VPan}  the other approach is presented for the  description of  generators of the field of $U$-invariants for the  adjoint representation.

	\section{The general construction of  $U$-projector}
	
	Let ${\mathfrak u}$ be an nilpotent Lie algebra over a field  $K$ of zero characteristic, ~~ $U=\exp(\ux)$ be the corresponding group,
	$\Ac$ a commutative associative finitely generated algebra defined over the field  $K$ and without zero devisors,
	~$\Fc$ be the field   $\mathrm{Frac}(\Ac)$. Let  $D$ be a homomorphism of the Lie algebra $\ux$ into the Lie algebra of locally nilpotent derivations of the algebra  $\Ac$. Then the group  $U$ acts on $A$ by the formula  $g(a)=\exp D_x(a)$,~ $g=\exp(x)$.  The ring (the field) of $U$-invariants coincides with the ring (the field)of  $\ux$-invariants.
	The field of $U$-invariants  $\Fc^U$ is a field of fractions of the ring of $U$-invariants  $\Ac^U$ ~\cite[Theorem 3.3]{VP}.
	It is known that the  field $\Fc^U$ is a pure transcendental  extension of the main field  $K$  (see \cite{M, Pu, V}).
	
	By definition, an $U$-\emph{projector} is an arbitrary homomorphism   $P:\Ac\to \Fc^U$ identical on $\Ac^U$.
	We are going to present the general construction of   $U$-projector. One can construct a free generator system  of generators of the field  $\Fc^U$ in terms of $U$-projector.
	
	Fix the chain of ideal   $\ux=\ux_n\supset\ux_{n-1}\supset\ldots\supset\ux_1\supset\ux_0={0}$, where
	$\mathrm{codim}(\ux_{i}, \ux_{i+1})=1$.  For each pair  $\ux_i\supset \ux_{i-1}$, the subalgebra of invariants  $\Ac^{\ux_i}$ is contained in  $\Ac^{\ux_{i-1}}$  $\Ac^{\ux_0}=\Ac$.
	Let ${i_1}$ be the least number such that $\Ac^{\ux_{i_1}}\ne \Ac$. Fix  $x_{i_1}\in\ux_{i_1}\setminus \ux_{{i_1}-1}$.
	\\
	\Lemma\Num\label{aaa}. There exist the elements  $a_{1,1}\in \Ac$,~ $a_{1,0}\in \Ac^{\ux}$, $a_{1,0}\ne 0$ obeying
	\begin{equation} \label{dxa} D_{x_{i_1}}(a_{1,1})= a_{1,0}.\end{equation}
	\Proof.  Since  $D_{x_{i_1}}$ is a local nilpotent derivation, there exist  $a_{1,1}\in \Ac$,,~ $a_{1,0}\in \Ac^{\ux_{i_1}}$,~ $a_{1,0}\ne 0$.
	It is sufficient to prove    $a_{1,0}\in \Ac^\ux$. Really,
	for any  $y\in\ux$  the element $[y,x_{i_1}]$  belongs to  $\ux_{{i_1}-1}$. Therefore,  $D_{[y,x_{i_1}]}(a)=0$ for all $a\in \Ac$.
	Then
	\begin{equation}\label{daaa}
	D_y (a_{1,0})= D_yD_{x_{i_1}}(a_{1,1})= D_{x_{i_1}}D_y(a_{1,1})+D_{[y,x_{i_1}]}(a_{1,1})= D_{x_{i_1}}D_y(a_{1,1}).\end{equation}
	The formula (\ref{daaa}) implies that the least containing  $a_{1,0}$ and $D_\ux$-invariant subspace $<a_{1,0}>$ is contained in  $\mathrm{Im} D_{x_{i_1}}$. Since $\ux$  is a finite dimensional nilpotent Lie algebra and every derivation $\ux$  is locally nilpotent, the subspace  $<a_{1,0}>$ is finite dimensional.
	The representation  $D_\ux$  has a triangular form in  $<a_{1,0}>$. There exists a nonzero vector annihilated by  all $D_y$, ~ $y\in \ux$. We may assume that it is   $a_{1,0}$. $\Box$
	
	Let $a_{1,0}$ and  $a_{1,0}$ be as in lemma  \ref{aaa}. The element  $Q_1=a_{1,1}a_{1,0}^{-1}$ belongs to the localization $\Ac_1$ of the algebra  $\Ac$  with respect to the denominator system generated by   $a_{1,0}$, it obeys  $D_{x_{i_1}}(Q_1)=1$.
	Consider the linear mapping  $S_1: \Ac\to \Ac_1^{\ux_{i_1}}$ defined by
	\begin{equation}\label{SSS}
	S_1(a)= a-D_{x_{i_1}}(a)Q_1+D^2_{x_{i_1}}(a)\frac{Q_1^2}{2!}+\ldots+D^k_{x_{i_1}}(a)\frac{Q_1^k}{k!}+\ldots.
	\end{equation}
	The mapping  $S_1$ is an algebra  homomorphism identical on  $\Ac^{\ux_{i_1}}$.
	One can extend the representation  $D$ of the Lie algebra $\ux$ to  $\Ac_1$, and each derivation  $D_y$ remains locally nilpotent on  $\Ac_1$. The mapping  $S_1$ to an homomorphism
	$\Ac_1$ to $\Ac_1^{\ux_{i_1}}$ identical on  $\Ac_1^{\ux_{i_1}}$.

	The subalgebra  $\Ac_1^{\ux_{i_1}}$ is invariant with respect to all  $D_x,~ x\in \ux$. Substitute
	$\Ac$ by $\Ac_1^{\ux_{i_1}}$,  and proceed as above.   Let  $i_2$ be least number obeying   $\Ac_1^{\ux_{i_2}}\ne \Ac_1^{\ux_{i_1}}$. Choose  $x_2\in \ux_{i_2}\setminus\ux_{i_2-1}$. As above there exist the elements   $a_{2,1}\in \Ac_1^{\ux_{i_1}}$ and  $a_{2,0}\in \Ac_1^{\ux}$, ~ $a_{2,0}\ne 0$ such that
	$$D_{x_{i_2}}(a_{2,1})= a_{2,0}.$$
	Let  $\Ac_2$ stands for localization of the algebra  $\Ac$ with respect to the denominator system generated by  $a_{1,0}, a_{2,0}$.
	Consider the element  $Q_2=a_{2,1}a_{2,0}^{-1}$ of the algebra  $\Ac_2$.  Similarly (\ref{SSS}), we construct the homomorphism
	$ S_2: \Ac_1^{\ux_{i_1}}\to \Ac_2^{\ux_{i_2}}$.
	
	Proceeding further  we obtain the chain
	$n\geq i_m>\ldots>i_2>i_1\geq 1$, the systems of elements   $a_{k,0}\in \Ac^\ux$  and  $a_{k,1}\in \Ac$ obeying  (\ref{dxa}), and the mappings  $S_1, S_2,\ldots S_m$, where each
	$S_k$ is a homomorphism  $\Ac_{k-1}^{\ux_{k-1}}\to \Ac_k^{\ux_k}$ identical on  $\Ac_{k-1}^{\ux_k}$. Denote by  $\Ac_*$ the localization of the algebra  $\Ac$ with respect to the denominator system generated by  $a_{1,0}, a_{2,0}, \ldots, a_{m,0}$. Consider the mapping
	$$P=S_m\circ\cdots\circ S_2\circ S_1.$$
	From all above we conclude.\\
	\Theorem\Num\label{acac}. The mapping  $P$ is a homomorphism of the algebra  $\Ac$ into     $\Ac_*^\ux$ identical on  $\Ac^\ux$. That is  $P$ is a  $U$-projector.\\
	\Remark. Since  $\{a_{k,0}\}\subset \Ac^\ux$, the projector  $P$ can be extended to a projector   $\Ac_*\to \Ac_*^\ux$ identical on  $\Ac_*^\ux$.
	
	Let the group  $U$  rationally act on the irreducible affine algebraic variety  $X$ defined over the field  $K$.
	Then the group  $U$ naturally acts in the algebra  of regular functions  $K[X]$ be the formula  $$T_gf(x)=f(g^{-1}x).$$
	Any regular function on  $X$ is contained in some  finite dimensional invariant subspace \cite[Lemma 1.4]{VP}. The representation  $D=d_e T$ of the Lie algebra  $\ux$ is nilpotent in this subspace.
	Therefore, for every  $x\in \ux$  the operator    $D_x$  is a locally nilpotent derivation of the algebra  $\Ac=K[X]$.
	As above   there exist  rational functions   $a_{k,0}(x)$, ~ $a_{k,1}(x)$,~~$ 1\leq k\leq m,$. Notice that $\{a_{k,0}(x)\}$ are  $U$-invariant, moreover $a_{1,0}$ and  $a_{1,1}$ are regular, and
	$a_{k,0}(x)$, ~ $a_{k,1}(x)$  belong to the localization of the algebra of regular functions with respect to the denominator subset generated by  $a_{i,0}(x)$, ~ $1\leq i\leq k-1$.
	Consider the  $U$-invariant open subset   $X_0=\{x\in X_0:~~ a_{k,0}(x)\ne 0, ~ 1\leq k\leq m\}$ and its subset
	$$\Sx=\{x\in X_0:~ a_{k,1}(x)=0,~~1\leq k\leq m\}.$$
	There is the restriction mapping  $\Res:K[X_0]\to K[\Sx]$. \\
	\Theorem\Num\label{gen}. Assume that the system  $\{a_{k,1}:~ 1\leq k\leq m\} $ generate the defining ideal of the subset  $\Sx$ in the algebra   $K[X_0]$.  Let  $b_1,\ldots, b_s\in \Ac$ be the system of elements such that $$\Res(b_1),\ldots,\Res(b_s), \Res(a_{0,0})^{\pm 1},\ldots, \Res(a_{m,0})^{\pm 1}$$ generate the algebra  $K[\Sx]$. Then  $P(b_1),\ldots, P(b_s),a_{0,0}^{\pm 1},\ldots, a_{m,0}^{\pm 1}$ generate the algebra of invariants $K[X_0]^U$. In particular,  $P(b_1),\ldots, P(b_s),a_{0,0}^{\pm 1},\ldots, a_{m,0}^{\pm 1}$ generate the field of invariants $K(X)^U$.\\
	\Proof.  Let $\Res_U$ stands for the  restriction of $\Res$ on the subalgebra of
	$U$-invariants.  Since all  $\Res(Q_i)$ annihilate on  $\Sx$, we have  $$\Res=\Res_U\circ P.$$ By the assumption,  $$\Res_UP(b_1),\ldots,\Res_UP(b_s),\Res_U(a_{0,0})^{\pm 1},\ldots, \Res_U(a_{m,0})^{\pm 1} $$ generate the algebra  $K[\Sx]$. To conclude the proof it is sufficient to prove that $\Res_U$ is an isomorphism  of the algebra  $K[X_0]^U$ onto $K[\Sx]$.
	Actually, the image $\mathrm{Im}(\Res_U)$ coincides with $K[\Sx]$. Let us show that  $\mathrm{Ker}(\Res_U)=0$. If $f$ is an $U$-invariant and $\Res_U(f)=0$, then
	$\Res(f)=0$. Hence $f=\phi_1 a_{1,1} + \ldots + \phi_m a_{m,1}$ for some  $\phi_1,\ldots,\phi_m\in K[X_0]$ and
	\begin{equation}\label{fpf}
	f=P(f)= P(\phi_1) P(a_{1,1}) + \ldots + P(\phi_m) P(a_{m,1}).\end{equation}
	Let us prove that  $P(a_{k,1})=0$ for each  $1\leq k\leq m$.
	Since the function  $a_{k,1}$ is $\ux_{i_k}$-invariant,  $S_i(a_{k,1})=a_{k,1}$ for all  $1\leq i<k$.
	The formula  (\ref{SSS}) implies  $S_k(a_{k,1})=0$.
	Then $$S_kS_{k-1}\cdots S_1(a_{k,1})=S_k(a_{k,1})=0$$
	for all  $1\leq k\leq m$.
	Therefore   $P(a_{k,1})=0$ for all  $k$. By the formula  (\ref{fpf}), we conclude  $f=0$. $\Box$
	
	\section{$U$-projectors for the adjoint representations}
	
	The goal of this section is to present an exact construction of the  $U$-projector for the adjoint representation reductive  split group. Let  $G$  be a  connected reductive split group over  a field $K$ of zero characteristic,  ~$\gx$ be its Lie algebra, ~ ~$\De$ be a root system with respect to the Cartan subalgebra  $\hx$ (respectively,  $\Dp$ be a set of positive roots),  $$\ux =\sum_{\al\in\De^+}KE_\al$$  be the standard maximal nilpotent subalgebra in $\gx$.  Via the Killing form we identify   $ \gx$ with  $\gx^*$, and the algebra of regular functions $K[\gx]$  with the symmetric algebra   $\Sc(\gx)=K[\gx^*]$. We extend the adjoint representation  $\ad_x$ of the Lie algebra  $\gx$ to the representation  in $\Sc(\gx)$  by derivations  $D_x(a)=\{x,a\}$, where  $x\in\gx$, ~ $a\in \Sc(\gx)$ and $\{\cdot,\cdot\}$ is natural Poisson bracket in  $\Sc(\gx)$. If  $x=E_\al$, then we denote  $D_\al=D_{x_\al}$.
	
	Let  $\xi=\xi_1$ be one of maximal roots in $\Dp$.
	Consider the subset  $\De_2$ of the root system that consists of all  $\al\in\De$ obeying  $(\al,\xi)=0$. The subset  $\De_2$  is a root system for the reductive subalgebra  $\gx_2=\{x\in \gx:~ [x, E_{\xi}]=0\}$. The subalgebra
	$\gx_2$ contains the maximal nilpotent subalgebra  $\ux_2$ spanned by  $E_\al$, where $\al$ runs through the set of positive roots  $\Dp_2$ in $\De_2$.
	
	The set  $\Dp$splits into two subsets  $\Dp=\Gamma\cup \Dp_2$, where
	$\Gamma$ consists of all roots $(\al,\xi)>0$. The subset  $\Gamma$ contains  $\xi$; denote  $\Gamma^0=\Gamma\setminus \{\xi\}$. For each  $\al\in\Gamma^0$, there exists a unique
	$\al'\in \Gamma^0$ such that  $\al+\al'=\xi$ (see \cite{J}). The subalgebra  $\nx=\nx_1$  spanned by $\{E_\al,~~ \al\in\Gamma\}$ is isomorphic to the Heisenberg algebra,  and it is an  ideal in $\ux$.
	
	The element $E_{\xi}$  is annihilated by all derivations  $D_x$, ~ $ x\in \ux$.  The derivations  $D_x$ can be extended to the derivations of the localization  $\Sc'(\gx)$ of the algebra  $\Sc(\gx)$ with respect to the denominator subset generated by  $E_{\xi} $. The element
	\begin{equation}\label{Qxi}
	Q_{\xi}=-\frac{1}{2}H_{\xi}E_{\xi}^{-1}\in \Sc'(\gx)
	\end{equation} obeys the equality  $D_{\xi}(Q_{\xi})=1$. Following the formula  (\ref{SSS}) we construct the projector   $S_\xi$ of the algebra  $\Sc'(\gx)$ onto the subalgebra of  $D_\xi$-invariants.
	For each  $\al\in \Gamma_0$,  the element
	\begin{equation}\label{Qal}
	Q_{\al}=-\frac{1}{N_{\al,\al'}}E_{\al'}E_{\xi}^{-1}\in \Sc'(\gx)
	\end{equation}
	obeys  $D_\al(Q_{\al})=1$. As above, for each  $\al\in\Gamma_0$, we construct the projector  $S_\al$.
	The mapping
	\begin{equation}\label{PPP}
	P_1=\left(\prod_{\al\in \Gamma_0} S_\al\right)\circ S_{\xi},
	\end{equation}
	is a homomorphism of the algebra     $\Sc(\gx)$ to the subalgebra of invariants $\Sc'(\gx)^{\nx_1}$ identical on   $\Sc(\gx)^{\nx_1}$.
	One can extend the mapping   $P_1$ to the projector of $\Sc'(\gx)$ onto the subalgebra  $\Sc'(\gx)^{\nx_1}$.
	Notice that the product in the formula  (\ref{PPP}) does not depent on the ordering of factors.\\
	\Lemma\Num\label{nnn}. Let  $\mx$  be a Heisenberg algebra,~ $[\mx,\mx]=Kz$, ~$\Sc(\mx)_z$ be
	the localization of $\Sc(\mx)$  with respect to the denominator subset generated by  $z$. Let  $\mx_0$ stand for the complimentary subspace  for  $Kz$  in $\mx$.  Then if  is a derivation $D$ of  the algebra $\mx$ obeys  $D(z)=0$ and $D(\mx_0)\subseteq \mx_0$, then there exists a unique element  $b_D\in (\mx_0)^2 z^{-1}\in \Sc(\mx)_z$ such that $D(a)=\{b_D,a\}$ for any $a\in \mx$.
	\\
	\Proof. The proof is similar  \cite[Lemma 4.6.8]{Dix}. $\Box$

	The element  $x\in \gx_2$  provides the derivation  $D_x(a)=[x,a]$ of the Heisenberg algebra  $\nx$ with  $D_x(E_\xi)=0$ and $D$ preserves  the subspace  $$\nx_0=\mathrm{span}\{E_\al:~\al\in \Gamma_0\}. $$ By the lemma  (\ref{nnn}), there exists a unique  $b_x\in \nx_0^2 E_\xi^{-1}$ such that $D_x(a)=\{b_x,a\}$.  Then for any element   $\widetilde{x}=x-b_x$ from  $\Sc'(\gx)$, we obtain  $[\widetilde{x}, \nx]=0$.
	For any  $x,y\in \gx_2$, we have
	$$\{\widetilde{x},\widetilde{y}\} =\{x-b_x,y-b_y\} = \{x,y-b_y\} =\{x,y\}  - \{x,b_y\}.$$
	Since  $\{x,y\}=[x,y]\in \gx_2$ è $ \{x,b_y\} \in \nx_0^2 E_\xi^{-1}$,  the lemma  \ref{nnn} implies  $  \{x,b_y\}= b_{[x,y]}$.  Hence
	\begin{equation}\label{thetaa}
	\widetilde{[x,y]} = \{\widetilde{x}, \widetilde{y}\}
	\end{equation}
	for any  $x,y\in \gx_2$. The subset  $\widetilde{\gx_2} = \{\widetilde{x}:~ x\in\gx_2\}$ is a Lie algebra with respect to the Poisson bracket  in  $\Sc'(\gx)$, and  it is isomorphic to   $\gx_2$.
	Applying   (\ref{thetaa}),  we obtain
	\begin{equation}\label{DDD}
	\widetilde{D_x(y)} = \widetilde{[x,y]} = \{\widetilde{x}, \widetilde{y}\} =\{x-b_x,y-b_y\} =  D_x(\widetilde{y})
	\end{equation}
	for any  $x,~y\in \gx_2$.
	
	Further we proceed analogically to what have been done for  $\gx$;  choose  a maximal root $\xi_2$ in  $\Dp_2$,  then we have got the subset of positive roots  $\Gamma_2$, the subalgebras  $\nx_2$,~ $\gx_3$ in $\gx_2$.
	The  $Q_{\xi_2}$  is defined similarly as for $\xi$, substituting  $H_\xi$ for $\widetilde{H}_{\xi_2}$ and  $E_\xi$ for $\widetilde{E}_{\xi_2}$   in the formula  (\ref{Qxi}). The element  $Q_{\xi_2}$ belongs to the subalgebra $\Sc''(\gx)$ that is the  localization of $\Sc(\gx)$ with respect to the denominator system generated by  $E_{\xi_1}$ and $\widetilde{E}_{\xi_2}$. The both elements  $E_{\xi_1}$ and $\widetilde{E}_{\xi_2}$ are $\ux$-invariants.
	
	Applying   (\ref{DDD}), we have
	$D_{\xi_2}(Q_{\xi_2})=1$.  Following the formula  (\ref{SSS}), we define the operator  $S_{\xi_2}$.
	Easy to show that if  $a\in \Sc'(\gx)^{\nx_1}$, then $S_{\xi_2}(a)$ also belogs to
	$\Sc'(\gx)^{\nx_1}$.
	
	Likewise  $S_{\xi_2}$ we define  $S_\al$ for each  $\al\in(\Gamma_2)_0 =\Gamma_2\setminus \xi_2$. The mapping
	\begin{equation}\label{PPPtwo}
	P_2=\left(\prod_{\al\in (\Gamma_2)_0} S_\al\right)\circ S_{\xi_2},
	\end{equation}
	is a homomorphism of the algebra    $\Sc(\gx)$ onto the subalgebra of invariants $\Sc''(\gx)^{\nx_2}$, and it is identical in $\Sc(\gx)^{\nx_2}$.
	The operator  $P_2$ extends to a projector  $\Sc''(\gx)$ into the subalgebra $\Sc''(\gx)^{\nx_2}$.
	If $a\in \Sc''(\gx)^{\nx_1}$, then $P_2(a)$ is invariant with respect to all $D_x$, ~$x\in\nx_1\oplus\nx_2$. Therefore, the mapping  $P_2\circ P_1$ is a homomorphism  of the algebras  $$\Sc(\gx)\to\Sc''(\gx)^{\nx_1\oplus\nx_2},$$
	and it if identical on  ${\nx_1\oplus\nx_2}$-invariants in  $\Sc(\gx)$.
	
	Continuing the process, we  obtain the chain of positive roots  $\xi=\xi_1, \xi_2, \ldots, \xi_m$, that is referred to as a \textit{ Kostant cascade}.
	The maximal nilpotent subalgebra $\ux$ decomposes into the sum of Heisenberg subalgebras  $$\ux=\nx_1\oplus\nx_2\oplus\cdots\oplus \nx_{m}$$ with  $[\nx_i,\nx_j]\subset\nx_i$ for all  $i<j$.
	
	We define the chain of subalgebras  $\gx=\gx_1\supset \gx_2\supset\ldots\supset \gx_m$ with the systems of positive roots  $\De=\De_1\supset\De_2\supset\ldots\supset \De_m$; each  $\xi_i$ is a maximal root in $\De_i$.
	By induction method with respect to   $i$, we define the projectors  $P_1, ~ P_2 ,\ldots, P_m$ such that each product   $P_i\circ\cdots\circ P_1$ is a projector of the algebra  $\Sc(\gx)$ into the subalgebra of $\nx_1\oplus\cdots\oplus \nx_i$-invariants in the localization of the algebra   $\Sc(\gx)$ with respect to the denominator system generated by  $$\Xi=\{E_{\xi_1}, \widetilde{E}_{\xi_2}, \ldots, \widetilde{\widetilde{E}}_{\xi_i}\}.$$
	Take $\Xi=\Xi_m$. Let $\Sc(\gx)_*$ stand for the localization of  $\Sc(\gx)$ with respect to the denominator system generated by  $\Xi$.\\
	\Theorem\Num\label{atwo}. \\
	1) The mapping  $P=P_m\circ\cdots\circ P_1$ is a homomorphism of the algebra   $\Sc(\gx)$ into  $\Sc(\gx)_*^U$ identical on     $\Sc(\gx)^U$, i.e.  $P$ is an  $U$-projector of the adjoint representation of the group  $G$.\\
	2) Let  $\{H_i\}$ be a basis orthogonal complement  to the Kostant cascade  in $\hx$. Then the system of elements  $$\{P(E_{-\al)}:~ \al\in\Dp\}\cup \{P(H_i)\}\cup  \Xi$$ freely generates  $\Sc(\gx)_*^U$  (it also freely generates the field of  $U$-invariants).

	\section{$U$-projectors for an arbitrary representations}
	
	In this section, we present the general scheme of construction of the $U$-projector for an arbitrary finite demensional representation.  As in the previous section  $K$ is a field of characteristic zero, ~ $\gx$ is a reductive split Lie algebra over the field $K$, ~  $G$ is a connected reductive group with the Lie algebra $\gx$, ~ $B$ is a Borel subgroup in  $G$ that containes the Cartan subgroup  $H$ and the maximal unipotent subgroup $U=B'$, their Lie algebras are  $\bx$,~ $\hx$, ~$\ux$.
	
	Let $V$ be an arbitrary finite dimensional representation of the group $G$. This representation
	defines the representation  $f(v)\to f(g^{-1}v)$ in the algebra $\Ac=K[V]$. The corresponding representation of the Lie algebra  Ëè $\gx$ is realized in this space by the formula  $$ D_x f(v)=-f(xv), ~ x\in \gx.$$
	For any  $x\in \ux$, the operator  $D_x$  is a locally nilpotent derivation of the algebra $\Ac$.
	
	Decompose  $V$ into a direct sum $V=W_0\oplus W_1$, where $W_0$ is an irreducible representation, and  $W_1$ is its invariant complement. Suppose that  $\dim W_0>1$.
	Choose a lowest vector $v_0$ in  $W_0$. The stabilizer  $\px^-$ of the one dimensional subspace  $<v_0>$ is a parabolic subalgebra in  $\gx$ containing  $\hx$. Suppose that   $\px^-\ne \gx$. Acting by the Cartan involution $\theta$ on  $\px^-$,
	we obtain the parabolic subalgebra   $\px$. The intersection  $\gx_1=\px\cup\px^-$ is a Levi subalgebra in $\px$ (and in $\px^-$). Let $\mx=\mathrm{rad}(\px)$.
	
	Here and further,  $\gamma_1\leq \gamma_2$ is an ordering on the set of all weights such that  $\gamma_2-\gamma_1$  is a sum of simple roots with nonnegative coefficients.
	
	We extend   $v_0$  to the basis  $v_0,v_1,\ldots, v_k,\ldots, v_n$ in $V$ as follows \\
	1) ~if  $v_i\in W_0$ and   $v_j\in W_1$, then  $i<j$;\\
	2) ~ let  $E_{\al_1},\ldots, E_{\al_k}$ be the basis of $\mx$ over $K$, then $v_i=E_{\al_i}v_0$,~ $1\leq i\leq k$,  is a basis in  $\mx v_0$, and  $v_{k+1},\ldots, v_n$is a basis  in the  $\gx_1$-invariant complement for  $<v_0>\oplus \mx v_0$  in $V$;\\
	3)~ each vector   $v_i$ is a weight vector with respect to  $\hx$, moreover if $v_i, v_j\in W_0$ and $\wt(v_i)<\wt(v_j)$, then $i<j$ (i.e.  $\al_i<\al_j$ implies $i<j$).
	
	Let  $\omega_0,\omega_1,\ldots, \omega_k,\ldots, \omega_n$ be the dual basis.
	The linear form $\omega_0$ is invariant with respect to  $\ux$; we extend the operators  $D_x$,~ $x\in\ux$ to locally nilpotent derivations of the localization  $\Ac_{\omega_0}$ of the algebra  $\Ac$ with respect to the denominator system generated by  $\omega_0$.
	Our first goal is to construct a projector  $ \Ac\to \Ac_{\omega_0}^\mx$.
	
	Notice that for each  $i$ the linear form  $D_{E_i}(\omega_j)$ belongs to the subspace  $<\omega_0,\ldots,\omega_{j-1}>$. Moreover, if  $i>j$, then $D_{E_i}(\omega_j)=0$.
	In the case  $i=j$, we have  $D_{E_j}(\omega_j)=-\omega_0$. Then for the element
	$Q_j=-\omega_j\omega_0^{-1}$ and each $i>j$, we obtain
	
	\begin{equation}   D_{E_i}(Q_j)=\left\{ \begin{array}{ll} 0,& \mbox{if} ~ i>j,\\1, & \mbox{if} ~ i=j.\end{array}\right.
	\end{equation}
	
	For each $1\leq i\leq k$, we construct  the  mapping  $S_{\al_i}$ according to  (\ref{SSS}).
	Define the mapping
	\begin{equation}\label{Pzero}
	P_0=S_{\al_1}\circ\cdots\circ S_{\al_k}.
	\end{equation}
	\Lemma\Num.
	The mapping $P_0$ is a homomorphism   $\Ac$ to $\Ac_{\omega_0}^\mx$ identical on $\Ac^\mx$; its kernel is generated (as an ideal) by  $\omega_1,\ldots,\omega_k$.\\
	\Proof.  One can directly verify  $S_{\al_j}(\omega_j)=0$ for all  $\leq j\leq k$.
	Since  $D_{E_i}(\omega_j)=0$ for $i>j$, we have  $S_{\al_i}(\omega_j)=\omega_j$ and
	$$P_0(\omega_j) = S_{\al_1}\circ\cdots\circ S_{\al_k}(\omega_j)= S_{\al_1}\circ\cdots\circ S_{\al_j}(\omega_j)=0.$$
	As  $P_0$ is a homomorphism of algebras, the ideal  $I$ that is generated by  $\omega_1,\ldots,\omega_k$ belongs to the kernel of  $P_0$.  On the other hand, for $j=0$ or $j>k$, the image $P_0(\omega_j)$ can be written in the form  $P_0(\omega_j)=\omega_j+b$, where $b\in I$. Therefore  $\mathrm{Ker}(P_0) = I$.
	
	Since  $S_{\al_j}(a)=a$ for each $j$ and  $\mx$-invariant  $a$,  we have $P_0(a)=a$.
	Let us show that for any  $a\in \Ac_{\omega_0}$ the image $P_0(a)$ is  $\mx$-invariant. For each  $1\leq s \leq k$ denote  $P_0^{(s)}=S_{\al_s}\circ\cdots\circ S_{\al_k}$.
	We shall prove by induction on  $s$, beginning from  $s=k$ that $$D_{\al_s}(P_0^{(s)}(a))=\cdots = D_{\al_k}(P_0^{(s)}(a))=0.$$
	
	Indeed, for  $s=k$ we obtain  $D_{\al_k}(P_0^{(k)}(a))= D_{\al_k}(S^{(k)}(a))=0$.
	Suppose that the statement holds for  $s+1$; let us prove it for  $s$.
	Easy to see that  $$D_{\al_s}(P_0^{(s)}(a))= D_{\al_s}S_{\al_{s}}(P_0^{(s+1)}(a))=0.$$
	Let  $t>s$. By induction on $n$, one can easily prove that for elements of an arbitrary Lie algebra the following equality holds
	\begin{equation}\label{xyn}
	xy^n=(y-\ad_y)^n(x).
	\end{equation}
	Applying  (\ref{xyn}), we verify that there exist the operators  $L_1,\ldots, L_{k-t}$ obeying $$D_{\al_t} S_{\al_s}= S_{\al_s}D_{\al_t}+L_1D_{\al_{t+1}}+\ldots+L_{k-t}D_{\al_{k}}.$$
	Then
	$$ D_{\al_t}(P_0^{(s)}(a)) = S_{\al_s}D_{\al_t}(P_0^{(s+1)}(a))+
	L_1D_{\al_{t+1}}(P_0^{(s+1}(a))+\ldots+L_{k-t}D_{\al_{k}}(P_0^{(s+1}(a)).  $$
	According to the induction assumption  $D_{\al_t}(P_0^{(s)}(a))=0$. $\Box$
	
	Let $G_1$ be a subgroup in  $G$ which Lie algebra coincides with $\gx_1$. The projector  $P_0$ is invariant with respect to $G_1$. Indeed, since  $g_1\mx g_1^{-1}=\mx$ for any  $g_1\in G_1$, the projector  $g_1P_0g_1^{-1}$ has the same kernel and image as  $P_0$, and hence $g_1P_0g_1^{-1} = P_0$. This implies $g_1P_0(a) = P_0(g_1a)$.
	
	The group  $G_1$ acts on the space  $V_1 = <v_0, v_{k+1},\ldots, v_n>$. The algebra $\Ac_1=K[V_1]$  is the symmetric algebra $$\Sc(V_1^*)=K[\omega_0,\omega_{k+1},\ldots,\omega_n].$$ The homomorphism  $P_0$ is an isomorphism of  $K[\omega_0^{\pm 1},\omega_{k+1},\ldots,\omega_n]$ to the algebra   $\Ac_{\omega_0}^\mx$. Since  $P_0$ commutes with  $g_1$, the operator $P_0$ is  an isomorphism of the $G_1$-representations   $V_1^*$ and $P_0(V_1^*)$.
	
	Choose the lowest (for  $\gx_1$) vector $v_0'$ in $V_1$,  and continue the process  as above.  Finally, we obtain  the chain of subspaces  $V\supset V_1\supset\ldots \supset V_s$, and that of reductive subalgebras  $\gx\supset\gx_1\supset\ldots \supset\gx_s$, where  $\gx_s$-action in    $V_s$ is diagonalizable. We obtain the chain of lowest vectors
	$$v_0, v_0',\ldots v_0^{(s-1)}\subset V_s$$
	and the corresponding linear forms  $$f_0=\omega_0, f_1=\omega_0',\ldots, f_{s-1} = \omega_0^{(s-1)}\subset V_s^*\subset V^*.$$
	We extend  $\{f_1, \ldots, f_{s-1}\}$ to the basis   $\{f_1, \ldots, f_{s-1},\ldots, f_m\}$ in $V_s^*$.
	For each $1\leq i\leq s-1$, determine a homomorphism  $P_i$ defined on the localization  $\Ac_{\Lambda_i}$ of the algebra  $\Ac$ with respect to the denominator system generated by  $$\Lambda_i = \{f_0, P_0(f_1), \ldots, P_{i-1}\circ\cdots\circ P_0(f_{i-1})\}.$$
	Denote   $P =P_{s-1}$,~ and $\Ac_*=\Ac_{\Lambda_{s-1}}$.\\
	\Theorem\Num\label{athree}. \\
	1) The mapping   $P$ is a homomorphism of the algebra  $\Ac$ to $\Ac_*^U$ identical on  $\Ac^U$, i.e.  $P$ is a   $U$-projector for the representation of $G$ in $V$. \\
	2) The system of elements  $\Lambda \cup \{P(f_i): 1\leq i\leq s-1 \} $ freely generate  $\Ac_*^U$ (and it freely generate the field of  $U$-invariants).

	\section{U-projector on reductive group}
	
	Let  $G$ be as above (i.e. it is a  connected reductive split group defined over a field $K$ of zero characteristic). We consider the representation of the group $G$ in the space  $K[G]$ defined by the formula  $R_gf(s)=f(g^{-1} s g)$. Our goal is to construct the  $U$-projector for this representation in the algebra  $K[G]$.
	
	Let  $\gx$  be the Lie algebra of the group  $G$.
	Let $\Pi=\{\al_1,\ldots,\al_n\}$ be the system of simple roots in  $\Delta^+$, and $\Phi=\{\phi_1,\ldots,\phi_n\}$ be the system of fundamental weights, $\phi_i(\al_j)=\delta_{ij}$.
	In each fundamental representation  $V_{i}$ with the highest weight $\phi_{i}$, we choose the highest vector  $v_{i}^+$
	and the lowest vector $ v_{i}^-$.  In the dual space  $V_i^*$, we choose the highest and the lowest vectors  $l_{i}^+$  and $l_{i}^-$ such that  $(v_i^-,l_i^+) = (v_i^+, l_i^-)=1$. The matrix element  $d_i(g)=(gv_i^+,l_i^+)$, ~ $1\leq i\leq n$, is an $U$-invariant. Denote by $K[G]_*$ the localization of the algebra  $K[G]$ with respect to the denominator system generated by  $\{d_i(g):~ 1\leq i\leq n\}$.

	Order the positive roots $\Dp=\{\beta_1,\ldots,\beta_m\}$ such that  $\beta_t<\beta_s$ implies $t<s$. Each   $\beta_s\in \Dp$ is either simple (i.e.  $\beta_s$ coincides with some  $\al_{\nu(s)}\in \Pi$), or  $\beta_s=\al_{\nu(s)}+\beta_s'$ for some  $\al_{\nu(s)}\in \Pi$ and $\beta_s'\in \Dp$. For each $\beta_s$, let us fix $\al_{\nu(s)}\in\Pi$.
	Let us correspond to each  $\beta_s \in\Dp$ the matric element  $$ d_{\beta_s}(g)=(gv^+,E_{-\beta_s}l^+).$$
	Then for any $x\in\ux$, we obtain
	\begin{equation}\label{dxm}
	D_xd_{\beta_s}(g)=-(xgv^+,E_{-\beta_s}l^+)+(gxv^+,E_{-\beta_s}l^+)=(gv^+,xE_{-\beta_s}l^+).
	\end{equation}
	
	If $x=E_\beta$, then as above we simplify notation  $D_\beta = D_{E_\beta}$.
	The formula  (\ref{dxm}) implies  $D_{\beta_s} d_{\beta_t}(g) = 0$,  if $s>t$, and
	$$D_{\beta_s} d_{\beta_s}(g) = \phi(H_{\beta_s})d_{\nu(s)}(g). $$
	Then for $$Q_{\beta_s}(g)= d_{\beta_s}(g)(\phi (H_{\beta_s})d_{\nu(s)}(g))^{-1},$$  we have
	\begin{equation}\label{qdqd}D_{\beta_s} Q_{\beta_t}(g) = \left\{\begin{array}{ll} 0,~& \mbox{if} ~ s>t,\\
	1,~& \mbox{if} ~ s=t.\end{array}\right.\end{equation}
	For each  $1\leq s\leq m$ we construct the mapping   $S_{\beta_s}$ by the formula  (\ref{SSS}).   Define the operator
	$$P=S_{\beta_1}\circ\cdots\circ S_{\beta_m}.$$
	For each  $\beta_s\in\Dp$,  consider the matrix element  $$c_{\beta_s}=(gE_{-\beta_s}v^+,l^+). $$
	\Theorem\Num\label{afour}. \\
	1) The operator  $P$ is a homomorphism  $K[G]_*$ to  $K[G]_*^U$ identical on   $K[G]_*^U$,  i.e. it is an $U$-projector.\\
	2) The system of rational functions  $$\{d_i(g):~ 1\leq i\leq n\}~ \bigcup ~ \{P(c_{\beta_s})(g):~1\leq s\leq m\}$$ freely generate the algebra  $K[G]_*^U$ (it also freely generate the field   $K(G)^U$).\\
	\Proof.  The statement  1) follows from above. Let us prove  2).  It is sufficient to show that the selected  system of functions  satisfy the conditions of the theorem  \ref{gen}.\\
	1) The inequalities $d_i\ne 0$, ~ $1\leq i\leq n$, define the open Bruhat cell  $Bw_0B$.
	Let us show that the system  $\{d_{\beta_s}(g): ~1\leq s\leq m\}$ generate the defining ideal  $I_{w_0B}$ of the subset  $w_0B$ in $Bw_0B$.
	The element  $g$ belongs to  $Bw_0B$ if $g$ can be written in the form $g=aw_0bh$, where $a=\exp(x)\in U$,~  $b=\exp(y)\in U$, ~$h\in H$, and  $$x=\sum_{\al\in \Dp} x_\al E_\al,~~  y=\sum_{\al\in \Dp} y_\al E_\al.$$
	The ideal  $I_{w_0B}$ is generated by  $\{x_\beta:~ \beta\in\Dp\}$.
	
	On the other hand,  $$d_{\beta_s}(g)=(aw_0bhv^+, E_{-\beta_s}l^+) = \phi_{\nu(s)}(h)(v^-, a^{-1}E_{-\beta_s}l^+) = \phi_{\nu(s)}(h)f_s(x),$$
	where $f_s(x)$ is a polynomial in  $\{x_\al\}$ of the form  $$f_s(x)=cx_{\beta_s} +~ \mbox{polynomial ~ in~} \{x_{\beta_t}:~ t<s\}$$
	with $c\ne 0$.
	Therefore, the ideal  $I_{w_0B}$ is generated by  $\{d_{\beta_s}(g): ~1\leq s\leq m\}$.\\
	2) Let us show that  $K[w_0B]$ is generated by the restrictions of $\{d_i(g),c_{\beta_s}(b)\}$ on $w_0B$.  Indeed,
	$K[w_0B]$ is generated by  $\{\phi_i(h), y_{\beta_s}\}$.
	
	On the other hand,
	$d_i(w_0bh)=(w_0bhv^+, l^+)=\phi_i(h)$,
	$$
	c_{\beta_s}(w_0bh)=(w_0bhE_{-\beta_s}v^+,l^+)=\phi_{\nu(s)}(h)(bE_{-\beta_s}v^+,l^-)=\phi_{\nu(s)}f'_s(y),
	$$
	where $f'_s(y)$ in a polynomial in  $\{y_\al\}$ of the form $$f'_s(y)=cy_{\beta_s} +~ \mbox{polynomial ~ in~} \{y_{\beta_t}:~ t<s\}$$
	with $c\ne 0$.
	Hence $K[w_0B]$ is generated by the restrictions of  $\{d_i(g),c_{\beta_s}(b)\}$ on $w_0B$.
	$\Box$


\begin{thebibliography}{100}
		\bibitem{M}
		Miyata K. Invariants of certain groups, 
		Nagoya Math. J. - 1971. - Vol. 41. - P. 69-73.
		\bibitem{Pu}
		Pukanszky L, Le\c cons sur les repr\'{e}sentations des
		groupes.  Monographies de la Soci\'{e}t\'{e} Math\'{e}matique de France. - Paris: Dunod, 1967.
		\bibitem{VPan}
		Vyatkina K.A., Panov A.N. Field of U-invariants of adjoint representation of the group
		GL(n,K) // Matematicheskie Zametki. - 2013. - Vol.93. - P. 144 -147.
		\bibitem{V}
		E.B.Vinberg, Ratinality of the field of invariants of the triangular group,
		Vestnik of Moscow univ. Series 1. Mathematics, mechanics.  -  1982. - No. 2. - P.  23-24.
		\bibitem{VP}
		Vinberg E.B., Popov V.L. Invariant theory,
		\bibitem{J}
		Joseph A. A Preparation theorem for the prime spectrum of a semisimple Lie algebra// J. Algebra. - 1977.- Vol. 48. - P. 241-289.
		\bibitem{Dix}
		Dixmier J. Alg\'{e}bres enveloppantes. - Paris: Gauthier-Villars, 1974.
		
	\end{thebibliography}
\end{document}